\documentclass[twoside]{article} 

\usepackage{jwn-article}
\usepackage{fullpage}
\usepackage{amsfonts,amssymb}
\usepackage[T1]{fontenc}
\usepackage{concrete}
\usepackage{alltt}
\usepackage{hyperref}
\usepackage{stmaryrd}
\usepackage{theorem}

\newcommand{\pp}[1][]{\ensuremath{\mathrm{PP#1}}}
\newcommand{\lp}{\ensuremath{\mathrm{LP}}}
\newcommand{\lpl}[1][]{\ensuremath{\mathrm{LP#1}^-}}
\newcommand{\lpu}[1][]{\ensuremath{\mathrm{LP#1}^+}}
\newcommand{\lin}[1]{\ensuremath{\left\llbracket{#1}\right\rrbracket}}
\newcommand{\llb}[1][\cdot]{\ensuremath{\llfloor{#1}\rrfloor}}
\newcommand{\lub}[1][\cdot]{\ensuremath{\llceil{#1}\rrceil}}
\newcommand{\errlb}[1][\cdot]{\ensuremath{E^-(#1)}}
\newcommand{\errub}[1][\cdot]{\ensuremath{E^+(#1)}}

\renewcommand{\Re}{\mathbb{R}}
\renewcommand{\le}{\leqslant}
\renewcommand{\ge}{\geqslant}

\theoremheaderfont{\scshape}
\theorembodyfont{\itshape}
\newtheorem{proposition}{Proposition}
\theorembodyfont{\rmfamily}
\newtheorem{example}{Example}

\title{Bounded Global Optimization for Polynomial Programming \\
using Binary Reformulation and Linearization}

\author{Joseph W. Norman \\ University of Michigan \\ 
\texttt{jwnorman@umich.edu}}
\date{May, 2012}

\begin{document}

\maketitle

\begin{abstract}
  This paper describes an approximate method for global optimization
  of polynomial programming problems with bounded variables.  The
  method uses a reformulation and linearization technique to transform
  the original polynomial optimization problem into a pair of mixed
  binary-linear programs.  The solutions to these two integer-linear
  reformulations provide upper and lower bounds on the global solution
  to the original polynomial program.  The tightness of these bounds,
  the error in approximating each polynomial expression, and the
  number of constraints that must be added in the process of
  reformulation all depend on the error tolerance specified by the
  user for each variable in the original polynomial program.  As these
  error tolerances approach zero the size of the reformulated programs
  increases and the calculated interval bounds converge to the true
  global solution.
\end{abstract}

\section{Introduction}

This paper describes a method for approximating with interval bounds
the global optimum of a polynomial program with bounded real variables
and no special convexity properties.  Such problems arise in many
contexts including the analysis of parametric probability models
\cite{norman-thesis,norman-problogic}.  For that application and many
others it is sufficient to generate interval bounds on the global
minimum and maximum solutions of the polynomial optimization problems
involved; it is not required to compute those solutions with infinite
precision.  The successive approximation method presented here uses a
reformulation and linearization technique to create a hierarchy of
pairs of mixed binary linear programs whose solutions bound the true
global optimum of the original polynomial problem.

The method allows the user to adjust directly and incrementally how
many variables and constraints are added in the reformulation stage,
and to compute hard bounds on the error in approximating each
polynomial before solving the reformulated problem; thus the user can
control computation to trade the tightness of the generated interval
bounds against the time required to compute them.  The reformulated
problems can be solved by well-developed branch-and-bound or
branch-and-cut methods for mixed integer linear programming for which
codes are available in commercial and free optimization software.  The
general global polynomial optimization problem is NP-hard and it
remains an open question how well the bounds provided by the proposed
method converge in practice and how best to control computation to
produce the most useful answers with the smallest amount of
computation.

The main contributions of this method are control of computation and
definitive interval-based answers with predictable and explicit bounds
on potential errors.  Based on the parameters provided by the user,
the method computes bounds on the possible error in the linear
approximation of the objective and each constraint, before the
reformulated optimization problem is solved.  This predictability of
error bounds means that the method can generate hard interval bounds
on the global solution to each polynomial optimization problem,
conditioned on the feasibility of that problem.  Feasibility can
sometimes be confirmed or refuted by solving the reformulated
problems; and even when feasibility cannot be established
definitively, the method offers useful linear bounds on the set of
points that are potentially feasible.

\subsection{Preliminaries}

Let us consider a polynomial program \pp{}:
\begin{equation}
  \label{eq:pp}
  \begin{array}{r@{\quad}l}
    \mbox{minimize} & f(x_1,\ldots,x_n) \\
    \mbox{subject to} &
    g_1(x_1,\ldots,x_n) \le 0 ,\;
    g_2(x_1,\ldots,x_n) \le 0 ,\; \ldots, \;
    g_q(x_1,\ldots,x_n) \le 0 \\
    \mbox{and} &
    \alpha_1 \le x_1 \le \beta_1, \;
    \alpha_2 \le x_2 \le \beta_2, \; \ldots, \;
    \alpha_n \le x_n \le \beta_n
  \end{array}
\end{equation}
in which the objective $f : \Re^n \rightarrow \Re$ is a polynomial
function of the variables $x_1,\ldots,x_n$ as is each constraint
$g_j$.  The lower bound on each real variable $x_i$ is $\alpha_i$ and
the upper bound is $\beta_i$.  Each variable $x_i$ may take positive,
negative, or zero values.  Fractional objectives can be accommodated
using the Charnes-Cooper transformation \cite{charnes}.  It will be
demonstrated how to generate from \pp{} a pair of mixed integer-linear
programs such that the global minimum of \pp{} is bounded by the
solution to each linear approximation.

Let us use $\mathbf{x}=(x_1,\ldots,x_n)$ to represent the tuple of
variables used in the polynomial program \pp.  Let us use $\Phi_x$ to
represent the set of values of $\mathbf{x}$ that satisfy the bound
constraints included in the polynomial program \pp:
\begin{equation}
\label{eq:bbox}
\Phi_x = \left\{ (x_1,\ldots,x_n) \;:\;
\alpha_1 \le x_1 \le \beta_1, \ldots, \alpha_n \le x_n \le \beta_n \right\}
\end{equation}
The bound constraints $\Phi$ and the polynomial constraints $g_j \ge
0$ define a semialgebraic set of feasible values for $\mathbf{x}$.
Thus we can rewrite the polynomial optimization problem \pp{} as
follows:
\begin{equation}
  \label{eq:ppvec}
  \begin{array}{r@{\quad}l}
    \mbox{minimize} & f(\mathbf{x}) \\
    \mbox{subject to} &
    g_1(\mathbf{x}) \le 0, \;
    g_2(\mathbf{x}) \le 0, \; \ldots, \;
    g_q(\mathbf{x}) \le 0 \\
    \mbox{and} &
    \mathbf{x} \in \Phi_x
  \end{array}
\end{equation}
Let us use $f_{\pp}(\mathbf{x}^*)$ to denote the true global minimum
solution to the problem \pp{}, where $\mathbf{x}^*$ identifies a
feasible (though not necessarily unique) point at which that minimum
occurs.  The goal is to develop a pair of mixed integer-linear
programs \lpl{} and \lpu{} using new variables $\mathbf{w}$ with the
property that their solutions bound the true global minimum:
\begin{equation}
  f_{\lpl}(\mathbf{w}^*) \quad \le \quad
  f_{\pp}(\mathbf{x}^*) \quad \le \quad
  f_{\lpu}(\mathbf{w}^*)
\end{equation}
It is also desirable for the linear approximations to provide
trustworthy information about the feasibility of the original
nonlinear program, reporting whether feasibility has been confirmed,
refuted, or not yet determined.  The reformulation and linearization
method presented below accomplishes these goals.

The tuple $\mathbf{w}$ of variables used in the reformulated programs
contains three classes of variables: `unit variables' denoted
$u_1,u_2,\ldots,u_\phi$ of which each binary variable $u_i \in
\{0,1\}$ is allowed to take the integer value zero or one; `remainder
variables' denoted $r_1,r_2,\ldots\,r_n$ of which each real variable
$0 \le r_i \le 1$; and `unit-product variables' denoted
$y_1,y_2,\ldots,y_\psi$ of which each real variable $0 \le y_i \le 1$.
The number of unit variables is $\phi$ and the number of unit-product
variables $\psi$.  Therefore the total number of variables in
$\mathbf{w}$ used in reformulated programs is $n+\phi+\psi$, where $n$
is the number of variables in the original program \pp.  Thus the
reformulated variables $\mathbf{w}$ are given by:
\begin{eqnarray}
  \label{eq:wvars}
  \mathbf{w} & = & (
  u_1,u_2,\ldots,u_\phi; 
  r_1,r_2,\ldots,r_n; 
  y_1,y_2,\ldots,y_\psi
  )
\end{eqnarray}
As part of the reformulation process some number $\rho$ of linear
constraints are constructed and added to the problem; these
constraints define a feasible set $\Phi_w$ of values for the
reformulated variables in $\mathbf{w}$.  Thus the problem is lifted
from the space $\Re^n$ to the space $\Re^{n+\phi+\psi}$ and linearized
in that space.  Analogous to the bounding box $\Phi_x$ that limits the
feasible values of $\mathbf{x} \in \Re^n$ a new bounding box $\Phi_w$
limits the feasible values of $\mathbf{w} \in \Re^{n+\phi+\psi}$.

In order to control the linearization and reformulation process the
user supplies parameters $\sigma_1,\sigma_2,\ldots,\sigma_n$ where
each $\sigma_i$ determines how the corresponding original variable
$x_i$ is to be approximated.  The tightness of the interval bounds
computed on the global solution to \pp{}, as well as the number $\phi$
of unit variables, the number $\psi$ of unit-product variables, and
the number $\rho$ of constraints that must be used in the
reformulation process, all depend on these user-supplied parameters.

\section{Prior Work}

\subsection{Reformulation and Linearization}

The method presented here is based on the principles set out by Li and
Chang \cite{li}, which in turn reference the reformulation
linearization technique of Sherali and Tuncbilek \cite{sherali}; it is
also a progression of the linearization and bounding techniques
described in the author's doctoral dissertation \cite{norman-thesis}.
The main extension relative to the work of Li and Chang is that here
there is a defined relationship between the solution to each linear
reformulation and the true solution to the original polynomial
program: the minimum solution of one linear program is less than or
equal to the true global minimum, which is in turn less than or equal
to the solution of the other linear program.  The original work
presented in \cite{li} offers the exact solution to a problem related
to the original polynomial problem, but with the exact qualifications
on that relationship left unspecified.  These desired solution
properties are produced by explicit calculation of the possible error
in approximating each polynomial expression through linearization
combined with reasoning about the inequalities involved.
Additionally, in this work a different method is used to linearize
products of two or more variables.  Finally, the number of variables
and constraints necessary to reformulate an optimization problem by
this method is calculated.  As in the original work no assumptions are
made about the convexity of the feasible set or objective function in
the polynomial problem; hence the technique presented here works for
maximization as well as minimization of objectives.  For simplicity of
presentation only minimization problems are described below.

The \emph{reformulation-linearization technique} (RLT) described by
Sherali and Tuncbilek \cite{sherali} presents an approach to solving
general polynomial programming problems which is quite similar.  In
each case the problem is lifted from its original $n$-dimensional
space $\Re^n$ to a higher-dimensional space via the addition of new
variables, and linear constraints are imposed in that
higher-dimensional space to provide approximations; the original RLT
provides only outer bounds (e.g.\ upper bounds on the global maximum),
whereas the present work provides inner bounds as well.
Mechanically, the present work uses two kinds of linearization
constraints (let us call them \emph{real-product-mean} and
\emph{binary-product-sum}) which are different from the
\emph{bound-factor}, \emph{constraint-factor}, and
\emph{convex-variable-bounding} constraints used in the traditional
RLT.  In Sherali and Tuncbilek's work the new variables in the lifted
space map one-to-one to the monomial terms that comprise the basis of
the original polynomial space; thus the number of new variables is
fixed.  Depending on which classes of linearization constraints are
added, the number of additional constraints varies; the user may
choose to add whole classes of constraints or select one at a time in
an ad-hoc way.  However in the present work the number of additional
variables can vary as well as the number of additional constraints.
Furthermore, these are both determined in a strict algorithmic way
according to a small set of parameters provided by the user.

One benefit of the approach taken here is that the user can make a
sequence of reformulations at different granularities or resolutions.
Moreover, the reformulation and subsequent linearization are separated
from the branch-and-bound search used to solve the resulting linear
programs (which are mixed integer linear programs in my case); that
is, the branch/bound requirements are encoded in integer variables.
This allows the current method to take advantage of well-developed
codes for solving mixed integer linear programs.  This also allows the
user to control the size of each reformulated program and thus the
time required to solve it.  We shall present in more detail the
relation to the traditional reformulation linearization technique.

\subsection{Semidefinite Programming and Real Algebraic Geometry}

Besides the linear relaxations of the RLT, another option is to use
semidefinite program relaxations of the original polynomial program,
to compute a lower bound on the true minimum.  As described by
Lasserre \cite{lasserre-moments,lasserre} this approach is in fact a
generalization of the linear relaxation approach and is related to it
through the theory of moments and its dual theory of representation of
polynomials which have positive values over a semialgebraic set.
Additionally, Floudas and Visweswaran \cite{floudas-1,floudas-2}
present a technique to transform bilinear, quadratic, or polynomial
programs to a new problem with partitioned variables and a certain
convexity property; the transformed problem can be solved by a primal
and relaxed dual approach.

There are many methods to solve special cases of the polynomial
programming problem, e.g. those in which all polynomial functions have
positive values or other instances in which favorable convexity
conditions are met.  For the broader class of constrained global
optimization problems there are also a variety of methods, including
outer approximation and branch-and-bound techniques, as reviewed by
Horst and Tuy \cite{horst-tuy}.  In the classification of Horst and
Tuy, the present method is a successive approximation method.  It is
deterministic rather than stochastic.  The inner or outer bounds at
each stage are computed by relaxation of the original polynomial
program to a mixed binary linear program.  Branching and bounding can
be carried out to solve each of these linear programs; however the
branching and bounding occurs in the (binary part of) lifted space not
the original space of the polynomial problem.  The current method uses
only the primal form of the polynomial program and does not
incorporate the dual problem in its analysis.

Most successive approximation methods reformulate the problem in one
step into a convex (linear or semidefinite) form; here we use a
nonconvex intermediate (a mixed 0-1 integer linear program) which is
NP-hard to solve but for which a great deal of work has been done to
make algorithms that perform well on average.  Note that successive
approximation with semidefinite relaxations is one of the possible
approaches to solving mixed binary linear programs which is distinct
from the more common branch-and-bound or branch-and-cut methods.

It is not clear to me if one can recover from each SDP relaxation any
bound on the error in the approximation (perhaps using the primal and
dual solutions provides appropriate bounds on the global minimum), or
if it is possible to determine error bounds on the approximation
before solving it; it seems instead that after solving a particular
relaxation it may be verified that it is indeed a global solution.

\section{Binary Reformulation and Linearization}

Here we develop the method to create a pair of mixed binary linear
programs to bound the global solution to a polynomial program.

\subsection{Binary Expansion and Basic Inequalities}

\begin{proposition}[Reformulation with unit variables]
\label{prop:reformulate}
Consider a real-valued variable $x_i$ bounded by the constants
$\alpha_i \le x_i \le \beta_i$.  Given a positive constant $\kappa_i$
satisfying $0<\kappa_i\le{\beta_i-\alpha_i}$ it is possible to
represent $x_i$ as a sum that involves some number $\sigma_i \ge 0$ of
binary variables $u_{i,j} \in \{0,1\}$ and a nonnegative real variable
$r_i \in [0,1]$:
\begin{eqnarray}
x_i & = &
\alpha_i + \kappa_i \sum_{j=1}^{\sigma_i}{2^{j-1}u_{i,j}} +
\kappa_i r_i
\end{eqnarray}
Let us call each $u_{i,j}$ a \emph{unit variable} and $r_i$ the
\emph{remainder variable}; each $\kappa_i$ is the \emph{error limit}
on the corresponding variable $x_i$.  In general it is necessary to
add the constraint that the sum satisfies the original upper bound
$\beta_i$:
\begin{equation}
\label{eq:reform-bound}
\alpha_i + \kappa_i \sum_{j=1}^{\sigma_i}{2^{j-1}u_{i,j}} +
\kappa_i r_i
\quad \le \quad \beta_i
\end{equation}
Otherwise it may be possible to choose values of $\kappa_i$, $u_{i,j}$
and $r_i$ that would produce a value greater than $\beta$.
\end{proposition}
This first proposition is a restatement of Equation~2.1 in \cite{li},
with additional detail provided here.  The number $\sigma_i$ of unit
variables required to reformulate a variable $x_i$ in this way is
related to the error limit $\kappa_i$:
\begin{eqnarray}
  \label{eq:to-sigma}
  \sigma_i & = & \left\lceil \log_2 
  \left(
  \frac{\beta_i-\alpha_i}{\kappa_i} 
  + \delta_i
  \right)
  \right\rceil
\end{eqnarray}
where the brackets indicate rounding up and $\delta_i=1$ if $x_i$ is
discrete or $0$ if $x_i$ is continuous.  The lowest possible error
limit $\kappa^*_i$ for a given number $\sigma_i$ of unit variables is
given by:
\begin{eqnarray}
  \label{eq:to-kappa}
  \kappa^*_i & = & \frac{\beta_i-\alpha_i}{2^{\sigma_i}-\delta_i}
\end{eqnarray}
where again $\delta_i=1$ if $x_i$ is discrete.  For each variable
$x_i$ the user may specify the number $\sigma_i$ of unit variables and
have the system compute the corresponding error limit $\kappa^*_i$
using Equation~\ref{eq:to-kappa}; or the user may specify the desired
error limit $\kappa_i$ and have the system compute the number of unit
variables required by Equation~\ref{eq:to-sigma}.

Note that in the case that $\sigma_i=0$ the error limit
$\kappa_i=\beta_i-\alpha_i$ and no unit variables will be added; note
also that a discrete variable with $\sigma_i=0$ must have
$\alpha_i=\beta_i$ in which case $\kappa_i$ is not needed for
reformulation.  Furthermore the error limit $\kappa_i$ can be zero
only for a fixed variable $x_i$ with bounds $\alpha_i=\beta_i$.  Also,
it turns out that it is not necessary to add any unit variables for a
variable $x_i$ which appears only linearly in \pp; in this case we
will take $\sigma_i=0$ but $\kappa_i=\beta_i-\alpha_i$ (thus the
reformulation of such an $x_i$ is $\kappa_i r_i$ which is equal to
$(\beta_i-\alpha_i)r_i$ with $0 \le r_i \le 1$).  Finally, we can
model a discrete variable with steps $\kappa_i$ using the substitution
above but omitting the remainder variable $r_i$; such a variable would
take values $\{\alpha_i, \alpha_i+\kappa_i, \alpha_i+2\kappa_i,
\ldots, \alpha_i+\sigma_i\kappa_i\}$ with $\alpha_i+\sigma_i\kappa_i
\le \beta_i$.


\begin{example}
\label{ex:vars}
Consider the variables $\mathbf{x}=(x_1,x_2,x_3)$ bounded by $2 \le
x_1 \le 5$, $0 \le x_2 \le 10$, and $4 \le x_3 \le 8$.  Let us use
three unit variables to reformulate $x_1$ and two unit variables to
reformulate each of $x_2$ and $x_3$.  In other words $\sigma_1=3$,
$\sigma_2=2$, and $\sigma_3=2$.  According to
Equation~\ref{eq:to-kappa} the smallest possible error limits are then
$\kappa^*_1=0.375$, $\kappa^*_2=2.5$, and $\kappa^*_3=1$.  The
reformulated variables according to Proposition~\ref{prop:reformulate}
are:
\begin{equation}
\label{eq:ex1vars}
\begin{array}{rcl}
x_1 & = &
2 + 0.375 u_{1,1} + 0.75 u_{1,2} + 1.5 u_{1,3} + 0.375 r_{1}
\\ x_2 & = &
2.5 u_{2,1} + 5 u_{2,2} + 2.5 r_{2}
\\ x_3 & = &
4 + u_{3,1} + 2 u_{3,2} + r_{3}
\end{array}
\end{equation}
\end{example}
With the following constraints added by Equation~\ref{eq:reform-bound}:
\begin{equation}
\label{eq:ex1constr}
\begin{array}{rcl}
2 + 0.375 u_{1,1} + 0.75 u_{1,2} + 1.5 u_{1,3} + 0.375 r_{1}
& \le &
5
\\
2.5 u_{2,1} + 5 u_{2,2} + 2.5 r_{2}
& \le &
10
\\
4 + u_{3,1} + 2 u_{3,2} + r_{3}
& \le &
8
\end{array}
\end{equation}

\begin{proposition}[Linear bounds on real products]
\label{prop:continue}
Consider nonnegative real-valued variables $r_1,r_2,\ldots,r_n$ each
bounded by $0 \le r_i \le 1$; assume $n>0$.  The difference between
the mean and the product of these bounded variables is limited by the
number of variables:
\begin{equation}
0 
\quad \le \quad
\left( 
\frac{
r_1 + 
r_2 + \cdots + 
r_n 
}{n}
\right) 
- 
\left( r_1 r_2 \cdots r_n \right)
\quad \le \quad 
\left( \frac{n-1}{n} \right) 
\end{equation}
which can be rewritten as:
\begin{equation}
0 
\quad \le \quad
(r_1 + r_2 + \cdots + r_n) - n (r_1 r_2 \cdots r_n )
\quad \le \quad 
n-1
\end{equation}
\end{proposition}
This proposition is a generalization of Equation~2.5 in
\cite{li}\footnote{There is a typographical error in Equation~2.5 in
  \cite{li}: the final term should be $\frac{1}{2}\omega_1\omega_2$
  instead of $\frac{1}{4}\omega_1\omega_2$.} with the following proof
added here.  It is clear that the relationships are true when every
$r_i$ is zero, in which case both the mean and product are zero; and
when every $r_i$ is one, in which case both the mean and product are
one.  For values other than these the greatest difference between the
mean and the product will occur when one of the variables is zero and
the others are one; in that case the mean will be $(n-1)/n$ and the
product will be zero, satisfying the equations above.

\begin{proposition}[Linear bounds on binary products]
\label{prop:product}
Consider binary variables $u_1,u_2,\ldots,u_n$ with each $u_i \in
\{0,1\}$ and also one real-valued `remainder' variable $r$ that
satisfies $0 \le r \le 1$.  The following relationships hold between
the product of the binary variables, the remainder variable, and the
individual binary variables:
\begin{equation}
\begin{array}{rcl}
u_1 u_2 \cdots u_n r & \le & u_i, \; i = 1,\ldots,n 
\end{array}
\qquad
\begin{array}{rcl}
u_1 u_2 \cdots u_n r & \ge & r + u_1 + u_2 + \cdots + u_n - n \\
u_1 u_2 \cdots u_n r & \le & r
\end{array}
\end{equation}
If the real variable $r$ is omitted the appropriate relationships
between the sum and product of the unit variables are:
\begin{equation}
\begin{array}{rcl}
u_1 u_2 \cdots u_n & \le & u_i, \; i = 1,\ldots,n 
\end{array}
\qquad
\begin{array}{rcl}
u_1 u_2 \cdots u_n & \ge & 1 + u_1 + u_2 + \cdots + u_n - n \\
u_1 u_2 \cdots u_n & \le & 1
\end{array}
\end{equation}
\end{proposition}
This third proposition is a restatement of Proposition~1 from
\cite{li} and the proof that appears there.  It is clear that if any
unit variable $u_j=0$ then the product $u_1 u_2 \cdots u_n = 0$.  In
that case the sum $u_1+u_2+\cdots+u_n \le n-1$ as one of the terms is
zero and each of the others is not greater than one.  The constraints
above are satisfied in this case.  If every unit variable $u_j=1$ then
the sum $u_1+u_2+\cdots+u_n=n$ and the product $u_1 u_2 \cdots u_n=1$,
in which case the constraints are also satisfied.

Let us introduce for this discussion a \emph{unit-product} variable
$y=u_1 u_2 \cdots u_n r$ to represent the product of the binary
variables $u_i$ through $u_n$ and the real-valued variable $r$
described above; we can use $y=u_1 u_2 \cdots u_n$ for the product of
the binary variables alone when no real variable $r$ is included.
Each unit-product variable $y$ is continuous and bounded by $0 \le y
\le 1$.

\subsection{Reformulation and Linearization of Polynomials}

Now we have the tools to reformulate any polynomial expression from
its native $x$-variables into a linear function of the unit variables
$u$, the unit-product variables $y$, and the remainder variables $r$.
This procedure is a generalization of Propositions 2 and 3 from
\cite{li} which describe only the linearization of products of two or
three variables.

\subsubsection{General Polynomial Form}

Consider that any polynomial function $g(x_1,x_2,\ldots,x_n)$ can be
represented as the sum of several terms $c_k m_k$ where each term is
the product of a real coefficient $c_k$ and a monomial $m_k$, the
latter of which is the product of several $x$-variables:
\begin{eqnarray}
  \label{eq:polynomial}
  g(x_1,x_2,\ldots,x_n) & = & c_1 m_1 + c_2 m_2 + \cdots + c_t m_t
  ,\quad c_k \in \Re
  ,\quad m_k = \prod_{i \in I_k} x_i
\end{eqnarray}
Let us adopt the convention that the first monomial $m_1=1$ so that
the coefficient $c_1$ is the constant term in the polynomial $g$.
Each tuple $I_k$ of indices identifies the original $x$-variables
included in the product for the monomial $m_k$.  An index $i$ may
occur in $I_k$ more than once; if a tuple $I_k$ is empty then the
corresponding product of zero variables is taken to be unity.  Note
that the size $|I_k|$ of the tuple of indices for a monomial $m_k$ is
the degree of that monomial.  Let us use $\mathbf{m}$ to represent the
list of all monomials which occur in a particular optimization problem
\pp, and $t$ for the number of such monomials.

\begin{example}
\label{ex:pp1}
Consider the following polynomial program \pp[1] which is the example
problem $\mathrm{PP}(\Omega)$ presented in Sherali and Tuncbilek
\cite{sherali} and reproduced as Example~1 in Li and Chang \cite{li}:
\begin{displaymath}
\begin{array}{r@{\quad}l}
\mbox{minimize}: & \begin{array}[t]{@{}l@{}} 5 x_2 + x_3 + x_1^{2} - 2 x_1 x_2 - 3 x_1 x_3 + 5 x_2 x_3 \\ \mbox{}  - x_3^{2} + x_1 x_2 x_3 \end{array} \\
\mbox{subject to}:
& \begin{array}[t]{@{}l@{}} 4 x_1 + 3 x_2 + x_3 \end{array} \le 20 \\
& \begin{array}[t]{@{}l@{}} x_1 + 2 x_2 + x_3 \end{array} \ge 1 \\
\mbox{and}:
& 2 \le x_1 \le 5 \\
& 0 \le x_2 \le 10 \\
& 4 \le x_3 \le 8
\end{array}
\end{displaymath}
The list of monomial terms (basis) used in \pp[1] is:
\begin{eqnarray}
  \label{eq:pp1m}
  \mathbf{m} & = & \left(
1, x_1, x_2, x_1 x_2, x_3, x_1 x_2 x_3, x_1^{2}, x_1 x_3, x_2 x_3, x_3^{2}
  \right)
\end{eqnarray}
The first term $m_1=1$ has the empty tuple $I_{1}=()$ of indices and
degree $0$; the last term $m_{10}=x_3^2$ has the tuple
$I_{10}=(3,3)$ of indices and degree $2$.  The number of terms in
$\mathbf{m}$ is
\begin{math}
t =
10
\end{math}
\end{example}

\subsubsection{Reformulation of Monomials with Sums of Unit Variables}

Now consider just one monomial $m_k = \prod_{i \in I_k} x_i$ of a
polynomial $g$ as given in Equation~\ref{eq:polynomial}.  Using
Proposition~\ref{prop:reformulate} we can substitute a sum involving
several unit variables and a remainder variable for each original
variable $x_i$ included in the product that defines this monomial $m_k$
where each $\alpha_i$ is the lower bound of the corresponding variable
$x_i$ and each $\sigma_i$ is the number of unit variables necessary to
represent $x_i$ within a tolerance of $\kappa_i$.  Carrying out the
multiplication 
to distribute the product over the sum yields an expression for the
monomial $m_k$ of the form:
\begin{equation}
  \label{eq:zsum}
  m_k \quad = \quad
  \prod_{i \in I_k}{\left( 
    \alpha_i + \kappa_i \sum_{j=1}^{\sigma_i}{2^{j-1}u_{i,j}} +
    \kappa_i r_i
    \right)}
  \quad = \quad
  z_{k,1} + z_{k,2} + \cdots + z_{k,s_k}
\end{equation}
where each $z_{k,\ell}$ designates the $\ell$-th element in the sum that
constitutes the distributed product.  With reference to
Equation~\ref{eq:zsum} the number $s_k$ of such elements is limited by:
\begin{equation}
\label{eq:slimit}
\prod_{i \in I_k} (\sigma_i+1) 
\quad \le \quad s_k \quad \le \quad
\prod_{i \in I_k} (\sigma_i+2)
\end{equation}
as the sum used to replace each variable $x_i$ according to
Proposition~\ref{prop:reformulate} contains $\sigma_i$ unit variables,
one remainder variable, and one constant lower bound $\alpha_i$ if
that bound is nonzero.

Now, it is clear from the structure of the inner sum and outer product
shown in Equation~\ref{eq:zsum} that every element $z_{k,\ell}$ in the
sum must be the product of some numbers of constant lower bounds
$\alpha$, unit variables $u$, and remainder variables $r$:
\begin{eqnarray}
  \label{eq:zprod}
  z_{k,\ell} & = &
  \left( 
  \alpha_{i'_1}
  \alpha_{i'_2} \cdots
  \alpha_{i'_{n_\alpha}}
  \right)
  \left( 
  \kappa_{i_1} 2^{j_1-1}
  u_{i_1,j_1}
  \kappa_{i_2} 2^{j_2-1}
  u_{i_2,j_2} \cdots
  \kappa_{i_{n_u}} 2^{j_{n_u}-1}
  u_{i_{n_u},j_{n_u}}
  \right)
  \left( 
  \kappa_{h_1} r_{h_1}
  \kappa_{h_2} r_{h_2} \cdots
  \kappa_{h_{n_r}} r_{h_{n_r}}
  \right)
\end{eqnarray}
where $n_\alpha$ is the number of constants, $n_u$ the number of unit
variables, and $n_r$ the number of remainder variables.  Note with
reference to Equation~\ref{eq:zsum} that the number of each kind of
element is limited by the degree $d_k = |I_k|$ of the monomial $m_k$:
\begin{equation}
n_\alpha \le d_k ,\quad n_u \le d_k ,\quad n_r \le d_k
\end{equation}
Let us introduce the constant $a_{k,\ell}$ to simplify notation:
\begin{eqnarray}
\label{eq:akl}
a_{k,\ell} & = &
\frac{1}{n_r}
\left( \alpha_{i'_1} \alpha_{i'_2} \cdots \alpha_{i'_{n_\alpha}} \right)
\left( \kappa_{i_1} \kappa_{i_2} \cdots \kappa_{i_{n_u}} 
2^{j_1-1} 2^{j_2-1} \cdots 2^{j_{n_u}-1} \right)
\left( \kappa_{h_1} \kappa_{h_2} \cdots \kappa_{h_{n_r}}
\right)
\end{eqnarray}
where we take the product of an empty set of items to be unity and
where we omit division by $n_r$ if $n_r=0$.  With this constant
$a_{k,\ell}$ we can rewrite the element $z_{k,\ell}$ shown in
Equation~\ref{eq:zprod} as:
\begin{eqnarray}
\label{eq:term}
z_{k,\ell} & = & a_{k,\ell} n_r
\left( u_{i_1,j_1} u_{i_2,j_2} \cdots u_{i_{n_u},j_{n_u}} \right)
\left( r_{h_1} r_{h_2} \cdots r_{h_{n_r}} \right)
\end{eqnarray}
where $n_r$ is omitted if it is zero.  Note that $a_{k,\ell}$ may be
positive or negative.

\begin{example}
\label{ex:reformulate}
The ninth term in the list $\mathbf{m}$ of monomial terms used in the
problem \pp[1] shown in Equation~\ref{eq:pp1m} is $m_{9} = x_2 x_3$.
Substituting the reformulation for each variable given in
Equation~\ref{eq:ex1vars} this term becomes:
\begin{equation}
m_{9} =
\left(
\begin{array}[t]{@{}l@{}} 2.5 u_{2,1} + 5 u_{2,2} + 2.5 r_{2} \end{array}
\right)
\left(
\begin{array}[t]{@{}l@{}} 4 + u_{3,1} + 2 u_{3,2} + r_{3} \end{array}
\right)
\end{equation}
Carrying out the multiplication gives the following sum of
\begin{math}
12
\end{math}
terms:
\begin{eqnarray}
  \label{eq:m8}
  m_{9} & = &
  \begin{array}[t]{@{}l@{}}
    10 u_{2,1} + 2.5 u_{2,1}u_{3,1} + 5 u_{2,1}u_{3,2} + 2.5 u_{2,1}r_{3}
    + 20 u_{2,2} + 5 u_{2,2}u_{3,1} 
    \\ \mbox{} 
    + 10 u_{2,2}u_{3,2} + 5 u_{2,2}r_{3} 
    + 10 r_2 + 2.5 u_{3,1}r_2 + 5 u_{3,2}r_2 + 2.5 r_2 r_3
  \end{array}
\end{eqnarray}
These terms are the elements $z_{9,1}$ through $z_{9,12}$ as described
in Equation~\ref{eq:zsum}.  Let us examine a few of these elements.
For the fourth element $z_{9,4}=2.5 u_{2,1} u_{3,1}$ we have: the
number of unit variables $n_u=2$, the number of remainder variables
$n_r=0$, and the constant $a_{9,2}=2.5$.  For the sixth element
$z_{9,6}=2.5 u_{2,1}r_{3}$ we have: $n_u=1$, $n_r=1$, and
$a_{9,6}=2.5$.  For the last element $z_{9,12}=2.5 r_2 r_3$ we have:
$n_u=0$, $n_r=2$, and $a_{9,12}=\frac{1}{2}(2.5)=1.25$.
\end{example}

\subsubsection{Linearization of Products of Remainder and Unit Variables}

We can use the propositions above to linearize each element
$z_{k,\ell}$ of the sum that constitutes the reformulated monomial
term $m_k$ in the following way.  First let us introduce the notation
$\lin{z_{k,\ell}}$ with double square brackets to denote the
linearized version of $z_{k,\ell}$.  Let us begin to compute this
linearized form by replacing the product of remainder variables in
$z_{k,\ell}$ with their mean:
\begin{eqnarray}
  \lin{z_{k,l}} & = & a_{k,\ell} n_r
  \left( u_{i_1,j_1} u_{i_2,j_2} \cdots u_{i_{n_u},j_{n_u}} \right)
  \left( \frac{ r_{h_1} + r_{h_2} + \cdots + r_{h_{n_r}} }{n_r} \right)
\end{eqnarray}
Simplifying this expression yields:
\begin{eqnarray}
  \label{eq:term-linear}
  \lin{z_{k,l}} & = & a_{k,\ell}
  \left( u_{i_1,j_1} u_{i_2,j_2} \cdots u_{i_{n_u},j_{n_u}} \right)
  \left( r_{h_1} + r_{h_2} + \cdots + r_{h_{n_r}} \right)
\end{eqnarray}
Note that if the element $z_{k,\ell}$ contains no remainder variables
($n_r=0$) the linearized form is unchanged from the original.  Let us
now replace the product of the unit variables
\begin{math}
\left( u_{i_1,j_1} u_{i_2,j_2} \cdots u_{i_{n_u},j_{n_u}} \right)
\end{math}
and each remainder variable $r_h$ in Equation~\ref{eq:term-linear}
with a new unit-product variable $y_{k,\ell,h}$ such that:
\begin{eqnarray}
\label{eq:yklh}
y_{k,\ell,h} & = &
\left( 
u_{i_1,j_1}
u_{i_2,j_2} \cdots
u_{i_{n_u},j_{n_u}}
\right)
r_h
,\quad h=1,2,\ldots,n_r
\end{eqnarray}
The linearized element $\lin{z_{k,\ell}}$ can now be written as
weighted sum of several unit-product variables:
\begin{eqnarray}
\lin{z_{k,\ell}} & = &
a_{k,\ell}
\left(
y_{k,\ell,1} +
y_{k,\ell,2} + \cdots +
y_{k,\ell,n_r}
\right)
,\quad \mbox{if $n_r>0$}
\end{eqnarray}
In the case that there are no remainder variables in the element
$z_{k,\ell}$ and thus $n_r=0$ let us introduce a single
unit-product variable:
\begin{eqnarray}
\label{eq:ykl1}
y_{k,\ell,1} & = &
\left( 
u_{i_1,j_1}
u_{i_2,j_2} \cdots
u_{i_{n_u},j_{n_u}}
\right)
\end{eqnarray}
and express the linearized element $\lin{z_{k,\ell}}$ appropriately:
\begin{eqnarray}
\lin{z_{k,\ell}} & = & a_{k,\ell} y_{k,\ell,1}
,\quad \mbox{if $n_r=0$}	
\end{eqnarray}
To simplify notation let us introduce the variable $n_{k,\ell}$ which
equals $n_r$ if $n_r>0$, or one otherwise.  Thus each linearized
element $\lin{z_{k,\ell}}$ can be expressed:
\begin{eqnarray}
  \label{eq:zkl}
  \lin{z_{k,\ell}} & = &
  a_{k,\ell} \sum_{h=1}^{n_{k,\ell}} y_{k,\ell,h}
  ,\quad n_{k,\ell} = 
  \left\{
  \begin{array}{r}
    \mbox{$1$ if $n_r = 0$} \\
    \mbox{$n_r$ if $n_r > 0$}
  \end{array}
  \right.
\end{eqnarray}

\begin{example}
\label{ex:zkl-lin}
Let us return to the elements that comprise the expanded form of the
monomial term $m_9=x_2 x_3$ described in Example~\ref{ex:reformulate}
above.  For the element $z_{9,2}=2.5 u_{2,1} u_{3,1}$, the linearized
form $\lin{z_{9,2}}=2.5 y_{9,2,1}$ using the solitary unit-product
variable $y_{9,2,1}=u_{2,1} u_{3,1}$ introduced according to
Equation~\ref{eq:ykl1}.  For the element $z_{9,4}=2.5 u_{2,1}r_3$ the
linearized form $\lin{z_{9,4}}=2.5 y_{9,4,1}$ using the unit-product
variable $y_{9,4,1} = u_{2,1}r_3$.  For the element $z_{9,12}=2.5 r_2
r_3$ the constant $a_{9,12}=1.25$ according to Equation~\ref{eq:akl}
and the linearized form $\lin{z_{9,12}}=1.25(r_2+r_3)$ according to
Equation~\ref{eq:term-linear}.  According to Equation~\ref{eq:yklh} we
could introduce two trivial unit-product variables $y_{9,12,1}=r_2$
and $y_{9,12,2}=r_3$ and express the linearized term
\begin{math}
  \lin{z_{9,12}} = 1.25 \left( y_{9,12,1} + y_{9,12,2} \right)
\end{math}
using the standard form given in Equation~\ref{eq:zkl}.
\end{example}

\subsubsection{Linear Constraints on Unit-Product Variables}
\label{sec:constraints}

For each element $z_{k,\ell}$ of the sum shown in
Equation~\ref{eq:zsum}, one unit-product variable $y_{k,\ell,h}$ has
been introduced for each of the $n_r$ remainder variables in
Equation~\ref{eq:term-linear} (or a single product variable
$y_{k,\ell,1}$ if $n_r=0$).  In order to satisfy
Proposition~\ref{prop:product} it is necessary to add the following
constraints for each unit-product variable $y_{k,\ell,h}$ thus
introduced:
\begin{equation}
\label{eq:constraints}
\begin{array}{rcl}
y_{k,\ell,h} & \le & u_{i_1,j_1} \\
y_{k,\ell,h} & \le & u_{i_2,j_2} \\
& \vdots & \\
y_{k,\ell,h} & \le & u_{i_{n_u},j_{n_u}}
\end{array}
\qquad
\begin{array}{rcl}
y_{k,\ell,h} & \ge &
r_h + u_{i_1,j_1} + u_{i_2,j_2} + \cdots + u_{i_{n_u},j_{n_u}} - n_u \\
y_{k,\ell,h} & \le & r_h
\end{array}
\end{equation}
where $r_h$ refers to the remainder variable associated with the
unit-product variable $y_{k,\ell,h}$ and with the substitution of
unity in the place of $r_h$ if $n_r=0$.  The number of added
constraints for each element $z_{k,\ell}$ is $n_r(n_u+1)$.  For a
unit-product variable $y_{k,\ell,h}=r_i$ which equals some remainder
variable (i.e. $n_r=1$ and $n+u=0$) or a unit-product variable
$y_{k,\ell,h}=u_{i,j}$ which equals some single unit variable
(i.e. $n_r=0$ and $n_u=1$), the constraints implied by
Equation~\ref{eq:constraints} are trivial and need not be added; it is
sufficient to note the identities.  In fact, in such cases the
original variables can be used and new unit-product variables need not
be introduced.

Let us use $\Phi_w$ to designate the set of values of $\mathbf{w}$ that
satisfy the linear constraints shown in Equation~\ref{eq:constraints}
for all unit-product variables introduced, as well as the upper bound
constraints as shown in Equation~\ref{eq:reform-bound} for all
reformulated original variables.  Let us use $\rho$ to designate the
number constraints used to specify the set $\Phi_w$.

\begin{example}
Let us use the elements discussed in Examples \ref{ex:reformulate} and
\ref{ex:zkl-lin} to illustrate the constraints on unit-product
variables.  The complete linearization of the reformulated monomial
term $m_9$ is given by:
\begin{eqnarray}
\lin{m_9} & = &
\begin{array}[t]{@{}l@{}} 10 u_{2,1} + 20 u_{2,2} + 11.25 r_{2} + 1.25 r_{3} + 2.5 y_{9,2,1} + 5 y_{9,3,1} \\ \mbox{}  + 2.5 y_{9,4,1} + 5 y_{9,6,1} + 10 y_{9,7,1} + 5 y_{9,8,1} + 2.5 y_{9,10,1} + 5 y_{9,11,1} \end{array}
\end{eqnarray}
Compare this with the reformulated but not yet linearized term shown
in Equation~\ref{eq:m8}; note that the corresponding terms appear in
different orders in the two equations.  The linearized form above uses
the following non-trivial unit-product variables:
\begin{equation}
\begin{array}{rcl}
y_{9,2,1} & = & u_{2,1} u_{3,1}
\\ y_{9,3,1} & = & u_{2,1} u_{3,2}
\\ y_{9,4,1} & = & u_{2,1} r_{3}
\\ y_{9,6,1} & = & u_{2,2} u_{3,1}
\\ y_{9,7,1} & = & u_{2,2} u_{3,2}
\\ y_{9,8,1} & = & u_{2,2} r_{3}
\\ y_{9,10,1} & = & u_{3,1} r_{2}
\\ y_{9,11,1} & = & u_{3,2} r_{2}
\end{array}
\end{equation}
These unit-product variables are accompanied by the constraints
described in Equation~\ref{eq:constraints}.  For example, for the
unit-product variable $y_{9,2,1}=u_{2,1}u_{3,1}$,
Equation~\ref{eq:constraints} implies the following constraints:
\begin{equation}
\label{eq:mc-8-2}
\begin{array}{rcl}
u_{2,1} \ge y_{9,2,1} \\
u_{3,1} \ge y_{9,2,1} \\
\begin{array}[t]{@{}l@{}} u_{2,1} + u_{3,1} \end{array} \le \begin{array}[t]{@{}l@{}} 1 + y_{9,2,1} \end{array}
\end{array}
\end{equation}
Similarly for the unit-product variable $y_{9,4,1} =
u_{2,1}r_3$ Equation~\ref{eq:constraints} implies:
\begin{equation}
\label{eq:mc-8-4}
\begin{array}{rcl}
u_{2,1} \ge y_{9,4,1} \\
\begin{array}[t]{@{}l@{}} u_{2,1} + r_{3} \end{array} \le \begin{array}[t]{@{}l@{}} 1 + y_{9,4,1} \end{array} \\
r_{3} \ge y_{9,4,1}
\end{array}
\end{equation}
For $z_{9,12}$ the constraints from Equation~\ref{eq:constraints} on
the unit-product variables $y_{9,12,1}=r_2$ and $y_{9,12,2}=r_3$ are
trivial and need not be added.
\end{example}

\subsubsection{Reformulation and Linearization of Entire Polynomials}

Recall from Equation~\ref{eq:zsum} that a polynomial term $m_k$
representing the product of several variables $x_i$ can be represented
as the sum of several elements $z_{k,\ell}$
\begin{eqnarray}
  m_k & = & z_{k,1} + z_{k,2} + \cdots + z_{k,s_k}
\end{eqnarray} 
where each element $z_{k,\ell}$ is the product of some number of
constants, unit variables, and remainder variables as shown in
Equation~\ref{eq:zprod}.  Using the reformulation technique above a
linear approximation $\lin{z_{k,\ell}}$ can be generated for
each element $z_{k,\ell}$ in this sum; adding these approximations
yields a linear approximation $\lin{m_k}$ for the polynomial term $m_k$:
\begin{eqnarray}
  \lin{m_k} & = &
  \lin{z_{k,1}} + \lin{z_{k,2}} + \cdots + \lin{z_{k,s_k}}
\end{eqnarray}
From the above and Equation~\ref{eq:zkl} it is clear that the linear
approximation for each monomial term $m_k$ is a weighted sum of
unit-product variables $y_{k,\ell,h}$:
\begin{eqnarray}
  \lin{m_k} & = &
  \sum_{\ell=1}^{s_k} \sum_{h=1}^{n_{k,\ell}} a_{k,\ell} y_{k,\ell,h}
\end{eqnarray}
where $s_k$ is the number of elements in the sum and each $n_{k,\ell}$
is the number of remainder variables included in the product that
defines each element, or $1$ if there are none.  The reformulation
also requires constraints involving the unit-product variables
$y_{k,\ell,h}$, the unit variables $u_{i,j}$, and the remainder
variables $r_i$ as shown in Equation~\ref{eq:constraints}.  It is
clear that this reformulation technique could be applied successively
to each monomial $c_k m_k$ in a polynomial $g$ as given in
Equation~\ref{eq:polynomial}, to yield a linear approximation
$\lin{g}$ of that polynomial:
\begin{equation}
  \label{eq:gsum}
  \lin{g} \quad = \quad
  c_1 \cdot \lin{m_1} + c_2 \cdot \lin{m_2} + \cdots + c_{t} \cdot \lin{m_{t}}
  \quad = \quad
  \sum_{k=1}^{t} \sum_{\ell=1}^{s_k} \sum_{h=1}^{n_{k,\ell}} 
  c_k a_{k,\ell} y_{k,\ell,h}
\end{equation}

\subsubsection{Bounds on Error from Linearization}
\label{sec:error}

In this section we shall consider the error introduced by the linear
reformulation process described above.  The error is a function of the
error limits $\kappa_1,\kappa_2,\ldots,\kappa_n$ set by the user for
the variables $x_i$ used in the polynomial expression $g$.  It follows
from Equations \ref{eq:term} and \ref{eq:term-linear} that for $n_r>1$
the difference between a linearized element $\lin{z_{k,\ell}}$ and its
true value $z_{k,\ell}$ is exactly:
\begin{eqnarray}
  \label{eq:zkl-err}
  \lin{z_{k,\ell}} - z_{k,\ell} 
  & = &
  a_{k,\ell}
  \left( 
  u_{i_1,j_1}
  u_{i_2,j_2} \cdots
  u_{i_{n_u},j_{n_u}}
  \right)
  \left( 
  (r_{h_1} + r_{h_2} + \cdots + r_{h_{n_r}})
  - n_r (r_{h_1} r_{h_2} \cdots r_{h_{n_r}})
  \right)
\end{eqnarray}
Note that the constant $a_{k,\ell}$ contains the product of several
error limits $\kappa_i$ as shown in Equation~\ref{eq:akl}.  In the
case that $n_r=0$ or $n_r=1$ no error is introduced; error is
introduced only by linearizing elements with $n_r \ge 2$ that contain
the products of two or more continuous remainder variables.

According to Proposition~\ref{prop:continue} the difference
$\lin{z_{k,\ell}}-z_{k,\ell}$ shown in
Equation~\ref{eq:zkl-err} is bounded by:
\begin{equation}
  0 \quad \le \quad \lin{z_{k,\ell}} - z_{k,\ell} \quad \le \quad
  a_{k,\ell}
  \left( 
  u_{i_1,j_1}
  u_{i_2,j_2} \cdots
  u_{i_{n_u},j_{n_u}}
  \right)
  (n_r-1)
  , \quad \mathrm{if} \ a_{k,\ell} > 0
\end{equation}
As the product 
\begin{math}
  \left( u_{i_1,j_1} u_{i_2,j_2} \cdots u_{i_{n_u},j_{n_u}} \right)
\end{math}
of the unit variables must be zero or one these error bounds are
simplified to:
\begin{equation}
  \label{eq:err-pos}
  0 \quad \le \quad \lin{z_{k,\ell}} - z_{k,\ell} \quad \le \quad
  a_{k,\ell} \left( n_r - 1 \right)
  , \quad \mathrm{if} \ a_{k,\ell} > 0,\ n_r > 0
\end{equation}
with the corresponding bounds for a negative coefficient
$a_{k,\ell}$:
\begin{equation}
  \label{eq:err-neg}
  a_{k,\ell} \left( n_r - 1 \right)
  \quad \le \quad \lin{z_{k,\ell}} - z_{k,\ell} \quad \le \quad 0
  , \quad \mathrm{if} \ a_{k,\ell} < 0,\ n_r > 0
\end{equation}
Using these relationships it is possible to create a pair of functions
to bound the error in the linear approximation of a polynomial.  

\begin{proposition}[Scalar bounds on error in linear approximation]
For a polynomial $g$ as given in Equation~\ref{eq:polynomial} let us
define the lower error bound $\errlb[g]$ to be the amount by which the
linearization $\lin{g}$ might underestimate the true value of $g$ and
the upper error bound $\errub[g]$ to be the amount by which the
linearization $\lin{g}$ might overestimate the true value of $g$:
\begin{equation}
  \errlb[g] \; \le \; \lin{g} - g \; \le \; \errub[g]
\end{equation}
With reference to the error bounds on elements $z_{k,\ell}$ given in
Equations \ref{eq:err-pos} and \ref{eq:err-neg}, the sum of
unit-product variables $y_{k,\ell,h}$ that defines each linearized
item $\lin{z_{k,\ell}}$ shown in Equation~\ref{eq:zkl}, and the sum of
unit-product variables that defines each linearized polynomial $[g]$
shown in Equation~\ref{eq:gsum}, these error bounds $\errlb[g]$ and
$\errub[g]$ can be computed as follows.  The lower error bound
$\errlb[g]$ is given by:
\begin{eqnarray}
  \label{eq:err-lb}
  \errlb[g] & = &
  \left\{ 
  \sum_{k=1}^{t} \sum_{\ell=1}^{s_k}
  c_k a_{k,\ell} \left( n_{k,\ell} - 1 \right)
  \; : \; c_k a_{k,\ell} < 0
  \right\}
\end{eqnarray}
where the sum includes an element only when the product of
coefficients $c_k a_{k,\ell}$ is negative.  The corresponding
expression for the upper error bound $\errub[g]$ is:
\begin{eqnarray}
  \label{eq:err-ub}
  \errub[g] & = &
  \left\{ 
  \sum_{k=1}^{t} \sum_{\ell=1}^{s_k}
  c_k a_{k,\ell} \left( n_{k,\ell} - 1 \right)
  \; : \; c_k a_{k,\ell} > 0
  \right\}
\end{eqnarray}
including elements only for positive products of coefficients $c_k
a_{k,\ell}$.  Note that the product
\begin{math}
  c_k a_{k,\ell} \left( n_{k,\ell} - 1 \right)
\end{math}
is nonzero only when $n_{k,\ell}>1$ (equivalently when $n_r>1$).
\end{proposition}

\begin{example}
Returning to the element $z_{9,12}$ given in
Example~\ref{ex:reformulate}, you can see from the reformulations for
the individual variables $x_2$ and $x_3$ given in
Example~\ref{ex:vars} that the element $z_{9,12} = 2.5 r_2 r_3$ is in
fact the product
\begin{math}
\left( \kappa_2 r_2 \right) \left( \kappa_3 r_3 \right)
\end{math}.
Thus according to Equation~\ref{eq:akl} the constant
$a_{9,12}=\frac{1}{2}\kappa_2\kappa_3$.  Equation~\ref{eq:err-pos}
says that the difference between the linearized element and its
original reformulation must be bounded by:
\begin{displaymath}
  0 \quad \le \quad \lin{z_{9,12}}-z_{9,12} \quad \le \quad
  \textstyle\frac{1}{2} \kappa_2 \kappa_3 (2-1)
\end{displaymath}
Substituting $\kappa_2=2.5$ and $\kappa_3=1$ gives 
\begin{math}
  0 \; \le \; \lin{z_{9,12}}-z_{9,12} \; \le \; 1.25
\end{math}.
Error bounds for complete polynomials can be computed using Equations
\ref{eq:err-lb} and \ref{eq:err-ub}.  For example the bounds on the
objective function:
\begin{eqnarray*}
f(x_1,x_2,x_3) & = &
\begin{array}[t]{@{}l@{}} 5 x_2 + x_3 + x_1^{2} - 2 x_1 x_2 - 3 x_1 x_3 + 5 x_2 x_3 \\ \mbox{}  - x_3^{2} + x_1 x_2 x_3 \end{array}
\end{eqnarray*}
of the problem \pp[1] in Example~\ref{ex:pp1} turn out to be
\begin{math}
\errlb[f] =
0
\end{math}
and
\begin{math}
\errub[f] =
1.25
\end{math}.
\end{example}

\subsubsection{Linear Bounds on Polynomials}

The reformulation and linearization procedure above, along with the
computed error bounds, allow us to compute linear expressions that
provide upper and lower bounds on any polynomial $g(\mathbf{x})$.
These bounds are valid for all values of $\mathbf{x}$ within the
feasible set $\Phi_x$.

\begin{proposition}[Linear bounds on polynomials]
Consider a polynomial $g$ which is a function of the variables
$\mathbf{x}=(x_1,x_2,\ldots,x_n)$.  Let $\Phi_x$ denote the set of values
of $\mathbf{x}$ that satisfy the bound constraints $\alpha_i \le x_i
\le \beta_i$ as shown in Equation~\ref{eq:bbox}.  Let $\lin{g}$ denote
the linear reformulation of $g$ according to the procedure given
above.  Let us introduce the notation $\llb[g]$ for the linear lower
bound on the polynomial $g$:
\begin{displaymath}
\llb[g] = \lin{g} - \errub[g]
\end{displaymath}
And similarly $\lub[g]$ for the linear upper bound:
\begin{displaymath}
\lub[g] = \lin{g} - \errlb[g]
\end{displaymath}
The construction above guarantees that the linear bounds are correct
for all feasible values of the original variables $\mathbf{x}$ and the
corresponding values of the reformulated variables $\mathbf{w}$:
\begin{eqnarray}
\llb[g](\mathbf{w}) & \le & g(\mathbf{x(w)})
, \quad \forall \mathbf{x(w)} \in \Phi_x \\
\lub[g](\mathbf{w}) & \ge & g(\mathbf{x(w)})
, \quad \forall \mathbf{x(w)} \in \Phi_x
\end{eqnarray}
\end{proposition}
where $\mathbf{x(w)}$ is the point in original variables corresponding
to the reformulated point $\mathbf{w}$.  The construction guarantees
that every point $\mathbf{x}$ in $\Phi_x$ has at least one
corresponding point $\mathbf{w}$ in $\Phi_w$, and that each point
$\mathbf{w}$ maps to a unique point $\mathbf{x(w)}$; however there may
be several feasible points $\mathbf{w(x)}$ to represent any given
$\mathbf{x}$.  Note that the reformulation $\lin{g}$ is a linear
function of the unit variables, remainder variables, and unit-product
variables which constitute the vector $\mathbf{w}$ shown in
Equation~\ref{eq:wvars}.  As each error bound $\errlb[g]$ and
$\errub[g]$ is a real number, the linear bounds $\llb[g]$ and
$\lub[g]$ defined above are therefore also linear functions of the
variables in $\mathbf{w}$.

\subsection{Pair of Bounding Mixed Binary Linear Programs}

We can use the linear bounds on polynomials described above to
generate more and less restrictive versions of the original polynomial
program \pp{} (inner and outer approximations).  As each reformulated
program will contain several binary unit variables as well as
continuous remainder variables and unit-product variables, it will be
a mixed integer linear program whose integer variables are binary.
For the optimistic case in which a lower bound on the global minimum
is desired, the linear lower bound $\llb[f]$ on the objective function
$f$ should be used in the optimistic reformulated program
$\lpl_\sigma$.  Also, each constraint $g_j \le 0$ in the polynomial
program \pp{} should be replaced in the linear program $\lpl_\sigma$
with its linear lower bound $\llb[g_j] \le 0$, which is less
restrictive.  For standardization each constraint $g_j \ge 0$ with the
inequality in the opposite direction should be replaced with the
equivalent constraint $-g_j \le 0$, and each equality constraint
$g_j=0$ replaced with the equivalent pair of constraints $g_j \le 0$
and $-g_j \le 0$.  Thus the looser problem $\lpl_\sigma$ is given by:
\begin{displaymath}
\label{eq:lpl}
\begin{array}{r@{\quad}l}
\mbox{minimize} & \llb[f](\mathbf{w}) \\
\mbox{subject to} &
\llb[g_1](\mathbf{w}) \le 0 ,\quad
\llb[g_2](\mathbf{w}) \le 0 ,\quad \ldots ,\quad
\llb[g_q](\mathbf{w}) \le 0 \\
\mbox{and} &
\mathbf{w} \in \Phi_w
\end{array}
\end{displaymath}
Let us use $\mathbf{w}^-$ to designate the point at which the minimum
solution to the optimistic program $\lpl_\sigma$
occurs, and $\llb[f](\mathbf{w}^-)$ to denote the value of the
linearized version of the objective function at that point.  The
linear reformulation $\lpl_\sigma$ can be
infeasible only if the original polynomial program \pp{} is
infeasible.  It can happen that the relaxed program
$\lpl_\sigma$ is feasible although the original
program \pp{} is not.

Similarly, for the pessimistic case in which an upper bound on the
global minimum is desired, the objective function and each constraint
in \pp{} should be replaced with its linear upper bound to produce the
tighter reformulated linear program $\lpu_\sigma$:
\begin{displaymath}
\label{eq:lpu}
\begin{array}{r@{\quad}l}
\mbox{minimize} & \lub[f](\mathbf{w}) \\
\mbox{subject to} &
\lub[g_1](\mathbf{w}) \le 0 ,\quad
\lub[g_2](\mathbf{w}) \le 0 ,\quad \ldots ,\quad
\lub[g_q](\mathbf{w}) \le 0 \\
\mbox{and} &
\mathbf{w} \in \Phi_w
\end{array}
\end{displaymath}
Let us use $\mathbf{w}^+$ to denote the point at which the minimum
solution $\lub[f](\mathbf{w}^+)$ of the pessimistic program
$\lpu_\sigma$ occurs.  The pessimistic
reformulation $\lpu_\sigma$ can be infeasible
even if the original program \pp{} is feasible; however if
$\lpu_\sigma$ is infeasible then \pp{} must be
infeasible as well.

We must use care in interpreting the solutions to the reformulated
linear programs $\lpl_\sigma$ and
$\lpu_\sigma$.  Let us use $\mathbf{x}^*$ to
designate a point at which the global minimum solution
$f(\mathbf{x}^*)$ to \pp{} occurs.  Let us say that a point
$\mathbf{w}$ in the reformulated variables is \emph{polynomial
feasible} if the corresponding point $\mathbf{x(w)}$ in the original
variables satisfies the bound constraints $\alpha_i \le x_i \le
\beta_i$ and the polynomial constraints $g_j \le 0$ in \pp{}.

The properties of the solutions to the reformulated programs are as
follows.  The solution to $\lpl_\sigma$ places a
lower bound on the true solution to \pp{} (there is no better
solution):
\begin{equation}
\label{eq:lp-llb}
\llb[f](\mathbf{w}^-) \le f(\mathbf{x}^*)
\end{equation}
If it happens that $\mathbf{w}^-$ is polynomial feasible, then the
value of the original polynomial objective function $f$ evaluated at
the corresponding point $\mathbf{x(w^-)}$ is an upper bound on the
true solution to \pp{} (there is a solution at least that good):
\begin{equation}
f(\mathbf{x}^*) \le f(\mathbf{x(w^-)})
, \quad \mbox{if} \;
\mathbf{x(w^-)} \in \Phi_x \; \mbox{and} \;
g_1(\mathbf{x(w^-)}) \le 0, \;
g_2(\mathbf{x(w^-)}) \le 0, \; \ldots, \;
g_q(\mathbf{x(w^-)}) \le 0
\end{equation}
However if the solution point $\mathbf{w}^-$ to
$\lpl_\sigma$ is not polynomial feasible then
the loose reformulation $\lpl_\sigma$ does not
provide any upper bound on the true solution to \pp{}; in that case it
is necessary to use an alternate means to generate an upper bound.
One way is to use the pessimistic reformulation
$\lpu_\sigma$ to compute an upper bound on the
global minimum solution to \pp.  If the tighter program
$\lpu_\sigma$ is feasible then its solution is
an upper bound on the global minimum:
\begin{equation}
f(\mathbf{x}^*) \le \lub[f](\mathbf{w}^+)
\end{equation}
If the reformulation \lpu{} is feasible, then by the construction of
\lpu{} the point $\mathbf{x(w^+)}$ in the original variables
corresponding to the solution point $\mathbf{w^+}$ to \lpu{} must be
polynomial feasible.  Therefore we can use the value of the original
objective function $f$ at that point as a tighter upper bound on the
global minimum:
\begin{equation}
f(\mathbf{x}^*) \le f(\mathbf{x(w^+)}) \le \lub[f](\mathbf{w}^+)
\end{equation}
In the case that \lpu{} is infeasible it does not provide an upper
bound on the global solution to \pp; and such a result does not prove
that \pp{} is infeasible.

\begin{example}
Let us consider problem \pp[1] from Example~\ref{ex:pp1}.  The
solution 
\begin{math}
\llb[f](\mathbf{w^-}) =
-124.799
\end{math}
to the optimistic reformulation 
\begin{math}
\lpl[1]_{(
3, 2, 2
)}
\end{math}
occurs at the point
\begin{math}
\mathbf{w^-}: (
u_{1,1} = 0, u_{1,2} = 1, u_{1,3} = 0, u_{2,1} = 0, u_{2,2} = 0, u_{3,1} = 1, u_{3,2} = 1, r_{1} = 0.666667, r_{2} = 0, r_{3} = 1
)
\end{math}.
This value $\llb[f](\mathbf{w^-})$ is a lower bound on the true global
minimum of \pp[1].  The corresponding point in the original variables
\begin{math}
\mathbf{x(w^-)}: (
x_1 = 3, x_2 = 0, x_3 = 8
)
\end{math}
happens to be polynomial feasible, satisfying the bound constraints
and two additional constraints in \pp[1].  The original objective
function evaluated at this point has value
\begin{math}
f(
x_1 = 3, x_2 = 0, x_3 = 8
) =
-119
\end{math}.
This value $f(\mathbf{x(w^-)})$ is an upper bound bound on the true global
minimum of \pp[1].  Thus the true global minimum solution to \pp[1]
lies in the interval
\begin{math}
[
-124.799, -119
]
\end{math}.
It turns out that the global minimum to \pp[1] is indeed
\begin{math}
-119
\end{math}
which we can prove by reformulating the problem using smaller error
limits and correspondingly more unit and unit-product variables.
\end{example}

\begin{example}
Let us consider a different polynomial optimization problem \pp[2]
which is used as Example~2 in Li and Chang \cite{li}:
\begin{displaymath}
\begin{array}{r@{\quad}l}
\mbox{minimize}: & \begin{array}[t]{@{}l@{}} 0.6224 x_3 x_4 + 19.84 x_1^{2} x_3 + 3.1661 x_1^{2} x_4 + 1.7781 x_2 x_3^{2} \end{array} \\
\mbox{subject to}:
& x_1 \ge 0.0193 x_3 \\
& x_2 \ge 0.00954 x_3 \\
& \begin{array}[t]{@{}l@{}} 1.33333 x_3^{3} \pi + x_3^{2} x_4 \pi \end{array} \ge 750.173 \\
& x_4 \le 240 \\
\mbox{and}:
& x_1 \in \{{1, 1.0625, 1.125, 1.1875, 1.25, 1.3125, 1.375}\} \\
& x_2 \in \{{0.625, 0.6875, 0.75, 0.8125, 0.875, 0.9375, 1}\} \\
& 47.5 \le x_3 \le 52.5 \\
& 90 \le x_4 \le 112 \\
& \pi = 3.14159
\end{array}
\end{displaymath}
Note that we can accommodate the discrete variables $x_1$ and $x_2$ by
reformulating each of them with binary unit variables as in
Proposition~\ref{prop:reformulate} but without the continuous
remainder variable $r_1$ or $r_2$.  The rest of the method works
without modification.  In this case the reformulation of each variable
is given by:
\begin{displaymath}
\begin{array}{rcl}
x_1 & = &
\begin{array}[t]{@{}l@{}} 1 + 0.0625 u_{1,1} + 0.125 u_{1,2} + 0.25 u_{1,3} \end{array}
\\ x_2 & = &
\begin{array}[t]{@{}l@{}} 0.625 + 0.0625 u_{2,1} + 0.125 u_{2,2} + 0.25 u_{2,3} \end{array}
\\ x_3 & = &
\begin{array}[t]{@{}l@{}} 47.5 + 0.00976562 u_{3,1} + 0.0195312 u_{3,2} + 0.0390625 u_{3,3} + 0.078125 u_{3,4} + 0.15625 u_{3,5} \\ \mbox{}  + 0.3125 u_{3,6} + 0.625 u_{3,7} + 1.25 u_{3,8} + 2.5 u_{3,9} + 0.00976562 r_{3} \end{array}
\\ x_4 & = &
\begin{array}[t]{@{}l@{}} 90 + 0.6875 u_{4,1} + 1.375 u_{4,2} + 2.75 u_{4,3} + 5.5 u_{4,4} + 11 u_{4,5} \\ \mbox{}  + 0.6875 r_{4} \end{array}
\end{array}
\end{displaymath}
using the error tolerances
\begin{math}
\left( \kappa_1,\kappa_2,\kappa_3,\kappa_4 \right) =
\left(
0.0625, 0.00976562, 0.6875, 0
\right)
\end{math}
and the corresponding numbers of unit variables
\begin{math}
\left( \sigma_1,\sigma_2,\sigma_3,\sigma_4 \right) =
\left(
3, 9, 5, 0
\right)
\end{math}
The loose reformulation 
\begin{math}
\lpl[2]_{(
3, 9, 5, 0
)}
\end{math}
of this problem \pp[2] has the
solution
\begin{math}
\llb[f](\mathbf{w^-}) =
6395.51
\end{math}
which is a lower bound on the true minimum $f(\mathbf{x}^*)$ as shown
in Equation~\ref{eq:lp-llb}.  This lower bound occurs at the
reformulated point
\begin{math}
\mathbf{w^-}: (
u_{1,1} = 0, u_{1,2} = 0, u_{1,3} = 0, u_{2,1} = 0, u_{2,2} = 0, u_{2,3} = 0, u_{3,1} = 0, u_{3,2} = 0, u_{3,3} = 0, u_{3,4} = 0, u_{3,5} = 0, u_{3,6} = 0, u_{3,7} = 0, u_{3,8} = 0, u_{3,9} = 0, u_{4,1} = 0, u_{4,2} = 0, u_{4,3} = 0, u_{4,4} = 0, u_{4,5} = 0, r_{1} = 0, r_{2} = 0, r_{3} = 0, r_{4} = 0, r_{5} = 0
)
\end{math}.
The corresponding point
\begin{math}
\mathbf{x(w^-)}: (
x_1 = 1, x_2 = 0.625, x_3 = 47.5, x_4 = 90, \pi = 3.14159
)
\end{math}
in the original variables happens to be polynomial feasible; the
original objective function evaluated at this point has
value
\begin{math}
f(
x_1 = 1, x_2 = 0.625, x_3 = 47.5, x_4 = 90, \pi = 3.14159
) =
6395.51
\end{math}.
Thus the true global minimum solution to \pp[2] lies in the interval
\begin{math}
[
6395.51, 6395.51
]
\end{math}.
This is a better solution than computed in Li and Chang \cite{li}.  More
importantly, the method presented here assures that this is in fact
the global minimum as there can exist no better solution than the
lower bound
\begin{math}
6395.51
\end{math}.
\end{example}

\begin{example}
We consider the problem \pp[3] which is Problem~338 in Schittkowski
\cite{schittkowski} and Example~3 in \cite{li}:
\begin{displaymath}
\begin{array}{r@{\quad}l}
\mbox{minimize}: & \begin{array}[t]{@{}l@{}} -x_1^{2} - x_2^{2} - x_3^{2} \end{array} \\
\mbox{subject to}:
& \begin{array}[t]{@{}l@{}} \frac{1}{2} x_1 + x_2 + x_3 \end{array} = 1 \\
& \begin{array}[t]{@{}l@{}} x_1^{2} + \frac{2}{3} x_2^{2} + \frac{1}{4} x_3^{2} \end{array} = 4 \\
\mbox{and}:
& -2 \le x_1 \le 2 \\
& -2.45 \le x_2 \le 2.45 \\
& -4 \le x_3 \le 4
\end{array}
\end{displaymath}
(The bounds on the variables were not specified in the original
formulation but they are implied by the equality constraints and made
explicit here.)  The program is reformulated using
\begin{math}
\left( \sigma_1,\sigma_2,\sigma_3 \right) =
\left(
7, 7, 7
\right)
\end{math}
and correspondingly
\begin{math}
\left( \kappa_1,\kappa_2,\kappa_3 \right) =
\left(
0.03125, 0.0382813, 0.0625
\right)
\end{math}.
The solution to the loose reformulation 
\begin{math}
\lpl[3]_{(
7, 7, 7
)}
\end{math}
is
\begin{math}
\llb[f](\mathbf{w^-}) =
-10.9965
\end{math}
which occurs at
\begin{math}
\mathbf{x(w^-)}: (
x_1 = -0.375, x_2 = -1.65897, x_3 = 2.84647
)
\end{math}.
However this point $\mathbf{x(w^-)}$ is not polynomial feasible; it
violates the second constraint in \pp[3].  It happens that the tight
reformulation 
\begin{math}
\lpu[3]_{(
7, 7, 7
)}
\end{math}
is infeasible, as is often the case for polynomial programs with
nonlinear equality (as opposed to inequality) constraints.  As an
alternate means of finding an upper bound on the true global minimum
to \pp[3] (and a polynomial-feasible point at which that upper bound
occurs) we can generate a focused problem that uses narrower ranges of
values for the variables $x_1$, $x_2$, and $x_3$ concentrated near the
polynomial-infeasible point $\mathbf{x(w^-)}$ where the solution to
\begin{math}
\lpl[3]_{(
7, 7, 7
)}
\end{math}
occurs.  If we choose a range focused on $\mathbf{x(w^-)} \pm
\kappa_i$ for each variable then the resulting polynomial program
is the following (let us call it \pp[3.2]):
\begin{displaymath}
\begin{array}{r@{\quad}l}
\mbox{minimize}: & \begin{array}[t]{@{}l@{}} -x_1^{2} - x_2^{2} - x_3^{2} \end{array} \\
\mbox{subject to}:
& \begin{array}[t]{@{}l@{}} \frac{1}{2} x_1 + x_2 + x_3 \end{array} = 1 \\
& \begin{array}[t]{@{}l@{}} x_1^{2} + \frac{2}{3} x_2^{2} + \frac{1}{4} x_3^{2} \end{array} = 4 \\
\mbox{and}:
& -0.40625 \le x_1 \le -0.34375 \\
& -1.69725 \le x_2 \le -1.62069 \\
& 2.78397 \le x_3 \le 2.90897
\end{array}
\end{displaymath}
This program \pp[3.2] can be reformulated using the same numbers of
unit variables
\begin{math}
\left( \sigma_1,\sigma_2,\sigma_3 \right) =
\left(
7, 7, 7
\right)
\end{math}
used for \pp[3] above, which now produce smaller error limits
\begin{math}
\left( \kappa_1,\kappa_2,\kappa_3 \right) =
\left(
0.000488281, 0.000598125, 0.000976562
\right)
\end{math}
due to the narrowed bounds on each variable.  The pessimistic
reformulated focused program
\begin{math}
\lpu[3.2]_{(
7, 7, 7
)}
\end{math}
has solution
\begin{math}
-10.9928
\end{math}
which is achieved at the polynomial-feasible point
\begin{math}
\mathbf{x(w^+)}: 
(
x_1 = -0.366211, x_2 = -1.6622, x_3 = 2.84531
)
\end{math}.
Combining the results of these reformulations 
\begin{math}
\lpl[3]_{(
7, 7, 7
)}
\end{math}
and 
\begin{math}
\lpu[3.2]_{(
7, 7, 7
)}
\end{math}
shows that the global solution to the original problem \pp[3] must lie
within the interval
\begin{math}
[
-10.9965
,
-10.9928
]
\end{math}.
Again, in contrast to the approximation method presented in Li and
Chang \cite{li}, the bounding approach presented here guarantees that
there cannot exist a better minimum than
\begin{math}
-10.9965
\end{math}.
\end{example}

\subsection{Alternative Formulation: Allowed Constraint Violation}

Another way to use the reformulation technique described above would
be to compute a single mixed binary linear program \lp{} from the
original program \pp{}, using the linearized version $\lin{g}$ of each
constraint $g$ and the linearized version $\lin{f}$ of the objective
$f$.  The error bounds on each $\lin{g}$ could then be used to
calculate the possible constraint violation $\tau_i$ for each
constraint $g_i \le 0$, and from these the maximum possible constraint
violation $\tau$ across all constraints could be computed.  The
interval $[z^-,z^+]$ would then contain the global optimum $z^*_\tau$
for the variant of the problem \pp{} in which each constraint is
\emph{nearly} satisfied (within the feasibility tolerance $\tau$).
The user could adjust the $\sigma_i$ parameters in the pre-solution
phase in order to achieve the desired feasibility tolerance.  Note
that in this alternative version, feasibility or infeasibility of the
reformulated optimization problem does not guarantee feasibility or
infeasibility of the original polynomial problem.

\section{Discussion}

\subsection{Problem Size}
\label{sec:size}

Let us now consider the number of variables and constraints that must
be added to a polynomial optimization problem \pp{} in the course of
reformulation and linearization as described above.  We assume that we
have a list $\mathbf{m} = \{m_1,m_2,\ldots,m_t\}$ of all the terms
used in monomials in the polynomial program \pp{} and we use $d$ to
denote the greatest degree of any term in $\mathbf{m}$.  We use
$\sigma = \sup_i{\sigma_i}$ to denote the largest number of unit
variables required to reformulate any variable $x_i$.

The original program \pp{} shown in Equation~\ref{eq:pp} contains:
\begin{itemize}
\item $n$ bounded real variables $x_1,x_2,\ldots,x_n$
\item $q$ polynomial constraints $g_1 \le 0, g_2 \le 0, \ldots g_q \le 0$.
\item One polynomial objective function $f$
\end{itemize}

In the reformulated program \lp{} there are:
\begin{itemize}
\item Exactly $n$ remainder variables $r_1,r_2,\ldots,r_n$ which are
real variables
\item Exactly $\sigma_1 + \sigma_2 + \cdots + \sigma_n$ unit
variables which are binary 0-1 variables:
\begin{displaymath}
\begin{array}{cccc}
u_{1,1} & u_{1,2} & \cdots & u_{1,\sigma_1} \\
u_{2,1} & u_{2,2} & \cdots & u_{2,\sigma_2} \\
& & \ddots & \\
u_{n,1} & u_{n,2} & \cdots & u_{n,\sigma_n}
\end{array}
\end{displaymath}
Let us use $\phi$ to represent the number of binary variables in a
reformulated mixed binary linear program \lp.  The text above shows
that $\phi = \sum_{i=1}^{n} \sigma_i$ which implies $\phi \le n
\sigma$.
\item Not more than $td(\sigma+2)^d$ unit-product variables which are
real variables.  For each monomial term $m_k \in M$:
\begin{displaymath}
\begin{array}{cccc}
y_{k,1,1} & y_{k,1,2} & \cdots & y_{k,1,n_{k,1}} \\
y_{k,2,1} & y_{k,2,2} & \cdots & y_{k,2,n_{k,2}} \\
& & \ddots & \\
y_{k,s_k,1} & y_{k,s_k,2} & \cdots & y_{k,s_k,n_{k,s_k}}
\end{array}
\end{displaymath}
where each $s_k$ is the number of terms
$z_{k,1},z_{k,2},\ldots,z_{k,s_k}$ needed to express the $k$-th
monomial term $m_k$ according to Equation~\ref{eq:zsum} and
$n_{k,\ell}$ is the number of remainder variables $n_r$ in the
$\ell$-th term $z_{k,\ell}$ or one if $n_r=0$.  Note that $n_r$ is
limited by the maximum degree $d$ of polynomial expressions in \pp:
$n_r{\le}d$.  Note also that the number of terms $s_k$ for each
monomial $m_k$ is limited by the degree $d_k$ of $m_k$ and the maximum
number $\sigma$ of binary variables used to represent any original
variable: $s_k \le (\sigma+2)^{d_k}$, using the product in
Equation~\ref{eq:slimit}.  Let us use $\psi$ to represent the number
of unit-product variables.  As shown here $\psi \le td(\sigma+2)^d$.
\item Not more than $td^2(d+1)(\sigma+2)^d$ product constraints, as
for each product variable $y_{k,\ell,h}$ added it is necessary to
add $n_r(n_u+1)$ constraints to satisfy
Proposition~\ref{prop:product} and both the numbers $n_r$ of
remainder variables (or one if there are none) and $n_u$ of unit
variables in any element $z_{k,\ell}$ are limited by the maximum
degree $d$ of polynomials.  Let us use $\rho$ to denote the number
of additional constraints.
\item Linearized versions of the original $q$ constraints.
\item The linearized objective function $\lin{f}$.
\end{itemize}
To summarize this using $O$-notation, the reformulated linear program
\lp{} will have a number of variables that is $O(td{\sigma^d})$
including $O(n\sigma)$ binary 0-1 variables; \lp{} will also have a
number of additional constraints (besides the original $q$ in \pp)
that is $O(td^3{\sigma^d})$.  Note that all of the added constraints
involve binary variables.

An important consideration is the time required to solve each mixed
binary linear program.  For a program with $\phi$ binary variables, in
the worst case a branch-and-bound algorithm would require enumeration
of the $2^\phi$ distinct combinations of $0$ and $1$ for each
variable, and solving a continuous linear program for each case.  As
$\phi$ is $O(n\sigma)$ the time required to solve each linear
reformulation is $O(2^{n\sigma})$ multiplied by the time required to
solve a standard linear program with $n+\phi+\psi$ variables and $\rho
\sim O(td^3\sigma^d)$ constraints.

\subsection{Conclusion}

This paper presented a reformulation and linearization technique to
generate an approximate solution to a polynomial optimization problem.
The approximate solution takes the form of interval bounds on the true
global optimum.  In the reformulation step each variable in the
original polynomial problem is replaced by a sum of binary and
continuous variables, the number of which depends on the error limit
specified by the user for each original variable.  In the
linearization step products of continuous variables are replaced by
sums of those variables, and constraints are added to the problem to
limit the differences between those sums and products.  Bounds on the
error introduced by linearization are computed, and with these bounds
a pair of mixed binary linear programs is created whose solutions
bound the solution to the original polynomial program.  The tightness
of the generated interval bounds depends on the error limits specified
by the user, which also determine the size of each reformulated
program and consequently the time required to solve it.

\begin{small}
  \raggedright
  \bibliographystyle{plain}
  \bibliography{../decisions}

  \vspace*{\baselineskip}
  \noindent\emph{Implementation note}:
\texttt{pqlsh-8.1.1+cplex+openmp 2012-05-29}
\end{small}

\end{document}